\documentclass[12pt,oneside,onecolumn,reqno]{amsart}
\usepackage[dvips]{epsfig}
\usepackage{amsgen, amstext,amsbsy,amsopn, amsthm, amsfonts,amssymb,amscd,amsmat
h,euscript,enumerate,url,verbatim,calc,oldgerm}
\usepackage{amsaddr}
\usepackage{bm}
\usepackage{latexsym}
\usepackage{color}
\usepackage{wrapfig}
\usepackage{graphicx}
\usepackage{epstopdf}
\usepackage{pifont}
\usepackage{setspace}
\usepackage{float}
\usepackage{filecontents}
\usepackage{blindtext}
\usepackage{csquotes}
\usepackage{cite}
\usepackage[caption = false]{subfig}

\theoremstyle{plain}
\newtheorem{theorem}{Theorem}[section]

\newtheorem{corollary}[theorem]{Corollary}

\theoremstyle{definition}

\newtheorem{example}[theorem]{Example}
\theoremstyle{remark}

\makeatletter

\newcommand{\Rmnum}[1]{\expandafter\@slowromancap\romannumeral #1@}
\makeatother

\begin{document}
\title{Bounds on Spreads of Matrices Related to Fourth Central Moment-II}
\author[R. Kumar]{R. Sharma$^{1}$, R. Kumar$^{2}$, \ R. Saini$^{3}$ \ and P.
Devi$^{4}$}
\address{$^{1,3, 4}$ Department of Mathematics, Himachal Pradesh University,
Shimla-171005, India \\
$^{2}$ Department of Mathematics, Dr B R Ambedkar National Institute of
Technology Jalandhar, Punjab-144011, India}
\email[R. Sharma]{rajesh.sharma.hpn@nic.in}
\email[R. Kumar]{ravithakur345@yahoo.com}

\begin{abstract}
We derive some inequalities involving first four central moments of discrete
and continuous distributions. Bounds for the eigenvalues and spread of a
matrix are obtained when all its eigenvalues are real. Likewise, we discuss
bounds for the roots and span of a polynomial equation.
\end{abstract}

\thanks{$^{\dag}$ AMS classification 60E15, 15A42, 12D10}
\keywords{Central moments, Trace, Positive linear functionals, Eigenvalues,
Roots, Span, Polynomial}
\maketitle



\section{Introduction}

This is the continuation of the work of Sharma et al. \cite{raj18}. It is
shown in \cite{raj18} that for both discrete and continuous random variable
in $[m,M]$, we have 
\begin{equation}
\begin{array}{lcl}
\label{ge1}\left\vert \mu _{3}\right\vert \leq \frac{\left( M-m\right) ^{3}}{%
6\sqrt{3}}\ \ \text{and}\ \ \mu _{4}\leq \frac{\left( M-m\right) ^{4}}{12}.
&  & 
\end{array}%
\end{equation}%
The inequalities \eqref{ge1} and the related Popoviciu inequality \cite%
{pop35}, 
\begin{equation}
\begin{array}{lcl}
\label{ge2}\mu _{2}\leq \frac{\left( M-m\right) ^{2}}{4} &  & 
\end{array}%
\end{equation}%
provide lower bounds for the range $r=M-m$ of the random variable in terms
of its central moments. \ These inequalities are also useful in many other
contexts. In literature, such inequalities are used to derive lower bounds
for the spread of a matrix, and span of a polynomial equation. The idea of
the spread is due to Mirsky \cite{mir56} and the notion of the span was
introduced by Robinson \cite{rob64}. For more detail and further related
topic see \cite{ara11, bha12, bha14, joh85, lek16, mer03, mir56, mir57,
raj10, raj12, raj13, raj15, raj16, raj17}.

\noindent In the present context we also need the following inequalities,
see \cite{bha00, mui66, raj10, raj12}, 
\begin{equation}
\begin{array}{lcl}
\label{age1}\mu _{2}\leq \left( M-\mu _{1}^{\prime }\right) \left( \mu
_{1}^{\prime }-m\right) &  & 
\end{array}%
\end{equation}%
\begin{equation}
\begin{array}{lcl}
\label{ge3}\frac{\mu _{2}^{2}-\left( \mu _{1}^{\prime }-m\right) ^{2}\mu _{2}%
}{\mu _{1}^{\prime }-m}\leq \mu _{3}\leq \frac{\left( M-\mu _{1}^{\prime
}\right) ^{2}\mu _{2}-\mu _{2}^{2}}{M-\mu _{1}^{\prime }} &  & 
\end{array}%
\end{equation}%
and 
\begin{equation}
\begin{array}{lcl}
\label{ge4}\mu _{2}+\left( \frac{\mu _{3}}{2\mu _{2}}\right) ^{2}\leq \frac{%
\left( M-m\right) ^{2}}{4}. &  & 
\end{array}%
\end{equation}%
A generalization of Samuelson's inequality \cite{sam68} due to Sharma and
Saini \cite{raj15} says that if $\overline{x}=\frac{1}{n}\sum_{i=1}^{n}x_{i}$
is the arithmetic mean and 
\begin{equation}
\begin{array}{lcl}
\label{age2}m_{r}=\frac{1}{n}\sum_{i=1}^{n}(x_{i}-\overline{x})^{r}\  &  & 
\end{array}%
\end{equation}%
is the $r$-th central moment of $n$ real numbers $x_{1},x_{2},...,x_{n}$,
then 
\begin{equation}
\begin{array}{lcl}
\label{ge5}m_{2r}\geq \frac{1+\left( n-1\right) ^{2r-1}}{n\left( n-1\right)
^{2r-1}}\left( x_{j}-\overline{x}\right) ^{2r} &  & 
\end{array}%
\end{equation}%
for all $j=1,2,...,n$ and $r=1,2,...$. For $r=1$, the inequality \eqref{ge5}
corresponds to Samuelson's inequality \cite{sam68}, and the complementary
inequalities due to Brunk \cite{bru59} assert that 
\begin{equation}
\begin{array}{lcl}
\label{mge27}m_{2}\leq \left( n-1\right) \left( M-\overline{x}\right) ^{2}%
\text{ \ and \ }m_{2}\leq \left( n-1\right) \left( \overline{x}-m\right)
^{2}. &  & 
\end{array}%
\end{equation}

\bigskip \noindent We here discuss some further extensions and applications
of the above inequalities to the field of theory of polynomial equations and
matrix analysis. Our first result gives an upper bound for $\mu _{4}$ in
terms of $m,M,\mu _{1}^{\prime }\ \text{and}\ \mu _{2}$ (Theorem 2.1, below)
and this provides further refinements of the inequality \eqref{ge2}
(Corollary 2.2-2.3). An extension of the second inequality \eqref{ge1} is
obtained for the fourth central moment $m_{4}$\ (Theorem 2.4). One more
lower bound for the range $r$ is obtained in terms of the second and fourth
central moment (Theorem 2.5). Some refinements of the first inequality %
\eqref{ge1} are given (Theorem 2.6). It is show that the inequalities %
\eqref{ge3} yield some more inequalities involving second and third central
moments (Theorem 2.7). The inequalities analogous to the inequality %
\eqref{mge27} involving second and fourth central moments are proved
(Theorem 2.8, Corollary 2.9-2.10). The upper (lower) bound for the smallest
(largest) eigenvalue of a complex $n\times n$ matrix are given when all its
eigenvalues are real as in case of Hermitian matrix (Theorem 3.1). The lower
bounds for the spread are obtained in terms of traces of $A^{r}$, $r=1,2,3,4$
(Theorem 3.3). Some results are extended for positive unital linear
functionals (Theorem 3.5). Likewise, we discuss bounds for the roots and
span of a polynomial equation (Theorem 4.1-4.2).

\section{Main Results}

\noindent It is enough to prove the following results for the case when $X$
is a discrete random variable taking finitely many values $%
x_{1},x_{2},...,x_{n}$ with probabilities $p_{1},p_{2},...,p_{n}\,$,
respectively. The arguments are similar for the case when $X$ is a
continuous random variable.

\begin{theorem}
\label{rt1} For $m\leq x_{i}\leq M,$ $i=1,2,...,n$, we have 
\begin{equation}
\begin{array}{lcl}
\label{mge1}\mu _{4}\leq \alpha \mu _{2}+\beta &  & 
\end{array}%
\end{equation}%
where 
\begin{equation}
\begin{array}{lcl}
\label{mge2}\alpha =\frac{3}{4}\left( M-m\right) ^{2}-2\left( \mu
_{1}^{\prime }-m\right) \left( M-\mu _{1}^{\prime }\right) &  & 
\end{array}%
\end{equation}%
and 
\begin{equation}
\begin{array}{lcl}
\label{mge3}\beta =\left( \mu _{1}^{\prime }-m\right) \left( M-\mu
_{1}^{\prime }\right) \left( \frac{1}{4}\left( M-m\right) ^{2}-\left( \mu
_{1}^{\prime }-m\right) \left( M-\mu _{1}^{\prime }\right) \right) . &  & 
\end{array}%
\end{equation}
\end{theorem}

\begin{proof}
Note that 
\begin{equation*}
\left( y+\frac{a+b}{2}\right) ^{2}\left( y-a\right) \left( y-b\right) \leq 0
\end{equation*}%
if and only if $a\leq y\leq b$. Applying this to $n$ numbers $y_{i}$'s, $%
a\leq y_{i}\leq b$; a little computation shows that 
\begin{equation}
\begin{array}{lcl}
\label{mge4}y_{i}^{4}\leq \left( \frac{3}{4}\left( a+b\right) ^{2}-ab\right)
y_{i}^{2}+\left( \frac{\left( a+b\right) ^{3}}{4}-\left( a+b\right)
ab\right) y_{i}-\frac{\left( a+b\right) ^{2}}{4}ab &  & 
\end{array}%
\end{equation}%
for all $i=1,2,...,n$. Multiplying both sides of \eqref{mge4} by $p_{i}\geq
0 $, and adding the $n$ resulting inequalities, we get that 
\begin{equation}
\begin{array}{lcl}
\label{mge5}\sum\limits_{i=1}^{n}p_{i}y_{i}^{4}\leq \left( \frac{3}{4}\left(
a+b\right) ^{2}-ab\right) \sum\limits_{i=1}^{n}p_{i}y_{i}^{2}+\left( \frac{%
\left( a+b\right) ^{3}}{4}-\left( a+b\right) ab\right)
\sum\limits_{i=1}^{n}p_{i}y_{i}-\frac{\left( a+b\right) ^{2}}{4}ab. &  & 
\end{array}%
\end{equation}%
For $y_{i}=x_{i}-\mu _{1}^{\prime },\ $we have $\sum%
\limits_{i=1}^{n}p_{i}y_{i}=0,$ $\sum\limits_{i=1}^{n}p_{i}y_{i}^{2}=\mu
_{2} $, $\sum\limits_{i=1}^{n}p_{i}y_{i}^{4}=\mu _{4}$, $a=m-\mu
_{1}^{\prime }\ $ and $b=M-\mu _{1}^{\prime }$. Using these facts and
choosing $y_{i}=x_{i}-\mu _{1}^{\prime }$ in \eqref{mge5}; the inequality %
\eqref{mge1} follows immediately.
\end{proof}

\bigskip \noindent Likewise, we can derive the upper bound for the fourth
moment about origin in terms of first and second moment, we have 
\begin{equation}
\begin{array}{lcl}
\label{mge6}\mu _{4}^{\prime }\leq \left( \frac{3}{4}\left( m+M\right)
^{2}-mM\right) \mu _{2}^{\prime }+\frac{1}{4}\left( m+M\right) \left(
M-m\right) ^{2}\mu _{1}^{\prime }-\frac{1}{4}mM\left( m+M\right) ^{2}. &  & 
\end{array}%
\end{equation}%
Further, the inequality \eqref{mge6} and hence the inequality \eqref{mge1}
becomes equality for $n=2$. In this case $\mu _{1}^{\prime }=mp_{1}+Mp_{2}$, 
$\mu _{2}^{\prime }=m^{2}p_{1}+M^{2}p_{2}$ and the right hand side
expression \eqref{mge6} equals $\mu _{4}^{\prime }=m^{4}p_{1}+M^{4}p_{2}$.

\begin{corollary}
\label{cr1} With the conditions as in Theorem 2.1, we have 
\begin{equation}
\begin{array}{lcl}
\label{mge7}\mu _{4}+3\mu _{2}^{2}\leq \left( M-m\right) ^{2}\left( \mu
_{1}^{\prime }-m\right) \left( M-\mu _{1}^{\prime }\right) \leq \frac{\left(
M-m\right) ^{4}}{4}. &  & 
\end{array}%
\end{equation}
\end{corollary}

\begin{proof}
The inequality \eqref{mge1} implies that 
\begin{equation}
\begin{array}{lcl}
\label{mge8}\mu _{4}+3\mu _{2}^{2}\leq \alpha \mu _{2}+\beta +3\mu _{2}^{2}.
&  & 
\end{array}%
\end{equation}%
From \eqref{age1} and \eqref{mge8}, we get that 
\begin{equation}
\begin{array}{lcl}
\label{mge9}\mu _{4}+3\mu _{2}^{2}\leq \alpha \left( M-\mu _{1}^{\prime
}\right) \left( \mu _{1}^{\prime }-m\right) +\beta +3\left( M-\mu
_{1}^{\prime }\right) ^{2}\left( \mu _{1}^{\prime }-m\right) ^{2}. &  & 
\end{array}%
\end{equation}%
Inserting the values of $\alpha $ and $\beta $ respectively from \eqref{mge2}
and \eqref{mge3} in \eqref{mge9}, and simplifying the resulting expression,
we immediately get the first inequality \eqref{mge7}. The second inequality %
\eqref{mge7} follows on using arithmetic mean-geometric mean inequality.
\end{proof}

\bigskip \noindent For $n=2$, the first inequality \eqref{mge7} becomes
equality. We have 
\begin{equation*}
\left( M-m\right) ^{2}\left( \mu _{1}^{\prime }-m\right) \left( M-\mu
_{1}^{\prime }\right) =p_{1}p_{2}\left( M-m\right) ^{4}=\mu _{4}+3\mu
_{2}^{2}.
\end{equation*}%
By \eqref{mge7}, we also get the following refinement of the inequality %
\eqref{ge2} for both discrete and continuos distributions,%
\begin{equation*}
\mu _{2}^{2}\leq \frac{\mu _{4}+3\mu _{2}^{2}}{4}\leq \frac{\left(
M-m\right) ^{4}}{16},
\end{equation*}%
cf. \cite{raj16}.

\noindent Let $m_{2}$ and $m_{4}$ respectively denote the second and fourth
central moment of $x_{1},x_{2},...,x_{n}$ as defined in \eqref{age2}. Then,
by \eqref{mge7},%
\begin{equation}
\begin{array}{lcl}
\label{mage5}m_{4}+3m_{2}^{2}\leq \frac{\left( M-m\right) ^{4}}{4}, &  & 
\end{array}%
\end{equation}%
and equality holds for $n=2$. The inequality \eqref{mage5} may be
strengthened when $n$ is odd.

\begin{corollary}
\label{cr2} For $m\leq x_{i}\leq M; i=1,2,...,n$, and with $n$ odd, the
inequality 
\begin{equation}
\begin{array}{lcl}
\label{mge10} m_{4}+3m_{2}^{2}\leq \frac{n^{2}-1}{4n^{2}}\left( M-m\right)
^{4} &  & 
\end{array}%
\end{equation}
holds true.
\end{corollary}

\begin{proof}
By the first inequality \eqref{mge7}, 
\begin{equation}
\begin{array}{lcl}
\label{mge11}m_{4}+3m_{2}^{2}\leq \left( M-m\right) ^{2}\left( M-\overline{x}%
\right) \left( \overline{x}-m\right) . &  & 
\end{array}%
\end{equation}%
Let $f\left( x\right) =\left( x-m\right) \left( M-x\right) $. Then the
function $f$ increases in $\left[ m,\frac{m+M}{2}\right] $, decreases in $%
\left[ \frac{m+M}{2},M\right] $ and attains its maximum at $x=\frac{m+M}{2}$%
. Thus, if $n$ is even the maximum value of $\left( \overline{x}-m\right)
\left( M-\overline{x}\right) $ is $\frac{\left( M-m\right) ^{2}}{4}$ and
this occurs when $x_{1}=x_{2}=...=x_{\frac{n}{2}}=m$ and $x_{\frac{n}{2}%
+1}=x_{\frac{n}{2}+2}=...=x_{n}=M$. But when $n$ is odd $\overline{x}\neq $ $%
\frac{m+M}{2}$ and $\left( \overline{x}-m\right) \left( M-\overline{x}%
\right) $ achieves its maximum when $x_{1}=x_{2}=...=x_{\frac{n-1}{2}}=m$
and $x_{\frac{n+1}{2}}=...=x_{n}=M$ or $x_{1}=x_{2}=...=x_{\frac{n+1}{2}}=m$
and $x_{\frac{n+3}{2}}=...=x_{n}=M$. Thus, the maximum is achieved either at 
\begin{equation*}
\overline{x}=\frac{\left( n-1\right) m+\left( n+1\right) M}{2n}\text{ \ or \ 
}\overline{x}=\frac{\left( n+1\right) m+\left( n-1\right) M}{2n},
\end{equation*}%
and in each case 
\begin{equation}
\begin{array}{lcl}
\label{mge12}\left( M-\overline{x}\right) \left( \overline{x}-m\right) \leq 
\frac{n^{2}-1}{4n^{2}}\left( M-m\right) ^{2}. &  & 
\end{array}%
\end{equation}%
The inequality \eqref{mge10} now follows from \eqref{mge11} and \eqref{mge12}%
.
\end{proof}

\bigskip \noindent In a similar vein we now prove an extension of the second
inequality \eqref{ge1} for the central moment $m_{4}$.

\begin{theorem}
\label{rt2} For $m\leq x_{i}\leq M,\ i=1,2,...,n$, and $I=\left\{
1,2,...,n-1\right\} $, we have 
\begin{equation}
\begin{array}{lcl}
\label{mge13}m_{4}\leq \max_{j\in I}\frac{j\left( n-j\right) \left(
n^{2}-3nj+3j^{2}\right) }{n^{4}}\left( M-m\right) ^{4}. &  & 
\end{array}%
\end{equation}
\end{theorem}

\begin{proof}
Let $f\left( x\right) =\left( x-\overline{x}\right) ^{4}$. The function $f$
is convex on the interval $\left[ m-\overline{x},\ M-\overline{x}\right] $.
Therefore, for $m\leq x\leq M$, 
\begin{equation*}
\left( x-\overline{x}\right) ^{4}\leq \frac{\left( M-\overline{x}\right)
^{4}-\left( m-\overline{x}\right) ^{4}}{M-m}\left( x-\overline{x}\right) -%
\frac{\left( m-\overline{x}\right) \left( M-\overline{x}\right) ^{4}-\left(
m-\overline{x}\right) ^{4}\left( M-\overline{x}\right) }{M-m},
\end{equation*}%
with equality if and only if $x=m$ or $x=M$. On applying this to $n$ numbers 
$x_{i}$'s, $m\leq x_{i}\leq M$ and using the arguments as in the proof of
Theorem 2.1, we obtain after simplification 
\begin{equation}
\begin{array}{lcl}
\label{mge14}m_{4}\leq \left( \overline{x}-m\right) \left( M-\overline{x}%
\right) \left( \left( \overline{x}-m\right) ^{2}+\left( M-\overline{x}%
\right) ^{2}-\left( \overline{x}-m\right) \left( M-\overline{x}\right)
\right) , &  & 
\end{array}%
\end{equation}%
with equality if and only if $x_{1}=x_{2}=...=x_{j}=m$ and $x_{j+1}=...=x_{n}=M$%
, $j=1,2,...,n-1.$ Hence, for some $j\in I$, the right hand side %
\eqref{mge14} achieves its maximum at 
\begin{equation*}
\overline{x}=\frac{jm+\left( n-j\right) M}{n}\text{.}
\end{equation*}%
For this value of $\overline{x}$ we have 
\begin{equation}
\begin{array}{lcl}
\label{mge15}\left( \overline{x}-m\right) \left( M-\overline{x}\right)
\left( \left( \overline{x}-m\right) ^{2}+\left( M-\overline{x}\right)
^{2}-\left( \overline{x}-m\right) \left( M-\overline{x}\right) \right) &  & 
\\ 
=\frac{j\left( n-j\right) \left( n^{2}-3nj+3j^{2}\right) }{n^{4}}\left(
M-m\right) ^{4}. &  & 
\end{array}%
\end{equation}%
The inequality \eqref{mge13} now follows from \eqref{mge14} and \eqref{mge15}%
.
\end{proof}

\begin{theorem}
\label{rt5} For $m\leq x_{i}\leq M,\ i=1,2,...,n$, we have 
\begin{equation}
\begin{array}{lcl}
\label{mage1}\mu _{2}\mu _{4}\leq \frac{4}{3^{5}}\left( M-m\right) ^{6}. & 
& 
\end{array}%
\end{equation}
\end{theorem}

\begin{proof}
We find from the inequality \eqref{mge1} that 
\begin{equation}
\begin{array}{lcl}
\label{mage2}\mu _{2}\mu _{4}\leq \left( \alpha \mu _{2}+\beta \right) \mu
_{2}. &  & 
\end{array}%
\end{equation}%
Combining \eqref{age1} and \eqref{mage2}, and inserting the values of $%
\alpha $ and $\beta $ respectively from \eqref{mge2} and \eqref{mge3}, we
get that 
\begin{equation}
\begin{array}{lcl}
\label{mage3}\mu _{2}\mu _{4}\leq \left( M-\mu _{1}^{\prime }\right)
^{2}\left( \mu _{1}^{\prime }-m\right) ^{2}\left( \left( M-m\right)
^{2}-3\left( \mu _{1}^{\prime }-m\right) \left( M-\mu _{1}^{\prime }\right)
\right) . &  & 
\end{array}%
\end{equation}%
The right hand side expression \eqref{mage3} achieves its maximum in the
interval $m\leq $ $\mu _{1}^{\prime }\leq M$ at $\mu _{1}^{\prime }=\frac{%
2m+M}{3}$ and $\mu _{1}^{\prime }=\frac{2M+m}{3},$ where its value is $\frac{%
4}{3^{5}}\left( M-m\right) ^{6}$. This proves the theorem.
\end{proof}

\bigskip \noindent The inequality \eqref{mage1} becomes equality for $n=2$
with $p_{1}=\frac{1}{3}$ and $p_{2}=\frac{2}{3}$ or $p_{1}=\frac{2}{3}$ and $%
p_{2}=\frac{1}{3}$. In this case, $\mu _{1}^{\prime }=\frac{2m+M}{3}$ or $%
\mu _{1}^{\prime }=\frac{2M+m}{3}$, $\mu _{2}=\frac{2}{9}(M-m)^{2}$ and $\mu
_{4}=\frac{2}{27}(M-m)^{4}$.\newline
\noindent It may be noted here that for $n=3$, $m_{4}=\frac{3}{2}m_{2}^{2}$
and the inequalities \eqref{mge10}, \eqref{mge13} and \eqref{mage1} become
same, $m_{4}\leq \frac{2}{27}\left( M-m\right) ^{4}$. \newline
\noindent We need Pearson's inequality \cite{pea16} in the proof of the
following theorem. This inequality gives a relation between skewness and
kurtosis of a distribution and can be written in the form%
\begin{equation}
\begin{array}{lcl}
\label{mage6}\frac{\mu _{3}^{2}}{\mu _{2}^{3}}+1\leq \frac{\mu _{4}}{\mu
_{2}^{2}}. &  & 
\end{array}%
\end{equation}

\begin{theorem}
\label{rt3} For $m\leq x_{i}\leq M,\ i=1,2,...,n$, \ the inequalities 
\begin{equation}
\begin{array}{lcl}
\label{mge16}\mu _{3}^{2}\leq \mu _{2}\mu _{4}-\mu _{2}^{3}\leq \frac{\left(
M-m\right) ^{6}}{108} &  & 
\end{array}%
\end{equation}%
and 
\begin{equation}
\begin{array}{lcl}
\label{mge17}\mu _{3}^{2}\leq \left( M-m\right) ^{2}\mu _{2}^{2}-4\mu
_{2}^{3}\leq \frac{\left( M-m\right) ^{6}}{108} &  & 
\end{array}%
\end{equation}%
hold true.
\end{theorem}

\begin{proof}
The first inequality \eqref{mge16} follows immediately from the inequality%
\eqref{mage6}.

\noindent Let $f\left( x\right) =ax-x^{3}.$ Then $f\left( x\right) $ has
maxima at $x=\sqrt{\frac{a}{3}}$ where its value is $\frac{2}{3\sqrt{3}}a^{%
\frac{3}{2}}.$ We thus have 
\begin{equation}
\begin{array}{lcl}
\label{mge18}\mu _{2}\mu _{4}-\mu _{2}^{3}\leq \frac{2}{3\sqrt{3}}\mu _{4}^{%
\frac{3}{2}}. &  & 
\end{array}%
\end{equation}%
The inequality \eqref{mge18} together with the second inequality \eqref{ge1}
yields the second inequality \eqref{mge16}.

\noindent The first inequality \eqref{mge17} follows immediately from the
inequality \eqref{ge4}. The second inequality \eqref{mge17} follows from the
fact that $f\left( x\right) =\left( M-m\right) ^{2}x^{4}-4x^{6}$ has maximum
at $x=\frac{M-m}{\sqrt{6}}$ and its maximum value is $\frac{\left(
M-m\right) ^{6}}{108}$.
\end{proof}

\bigskip \noindent For $n=2$, $\mu _{3}^{2}=\mu _{2}\mu _{4}-\mu
_{2}^{3}=p_{1}^{2}p_{2}^{2}(p_{1}-p_{2})^{2}(M-m)^{6}$ if and only if $p_{1}=%
\frac{1}{2}\pm \frac{1}{2\sqrt{3}}$ and $p_{2}=\frac{1}{2}\mp \frac{1}{2%
\sqrt{3}}$, and in this case inequalities \eqref{mge16} become equalities.
Likewise, the inequalities \eqref{mge17} reduce to equalities for $n=2$.

\noindent From \eqref{mge18} and the first inequality \eqref{mge16}, we find
that 
\begin{equation*}
\frac{\mu _{3}^{2}}{\mu _{2}^{3}}\leq \frac{2}{3\sqrt{3}}\frac{\mu _{4}}{\mu
_{2}^{2}}\sqrt{\frac{\mu _{4}}{\mu _{2}^{2}}}.
\end{equation*}%
This gives the following inequality involving skewness $\left( \alpha
_{3}\right) $ and kurtosis $\left( \alpha _{4}\right) $, 
\begin{equation*}
\frac{\alpha _{3}^{4}}{\alpha _{4}^{3}}\leq \frac{4}{27}.
\end{equation*}%
Further, from the first inequality \eqref{mge17}, we have 
\begin{equation*}
\alpha _{3}^{2}+4\leq q^{2},
\end{equation*}%
where $q$ is studentized range, $q=\frac{r}{\sqrt{\mu _{2}}}$.

\begin{theorem}
\label{rt4} For $m\leq x_{i}\leq M,\ i=1,2,...,n$, we have 
\begin{equation}
\begin{array}{lcl}
\label{mge19}-\frac{4}{27}\left( \mu _{1}^{\prime }-m\right) ^{5}\leq \mu
_{2}\mu _{3}\leq \frac{4}{27}\left( M-\mu _{1}^{\prime }\right) ^{5}, &  & 
\end{array}%
\end{equation}%
\begin{equation}
\begin{array}{lcl}
\label{mge20}\frac{1}{4}m^{2}\left( m+3\mu _{1}^{\prime }\right) \leq \mu
_{3}^{\prime }\leq \frac{1}{4}M^{2}\left( M+3\mu _{1}^{\prime }\right) , & 
& 
\end{array}%
\end{equation}%
\begin{equation}
\begin{array}{lcl}
\label{mage4}-\frac{1}{4}\left( \mu _{1}^{\prime }-m\right) ^{3}\leq \mu
_{3}\leq \frac{1}{4}\left( M-\mu _{1}^{\prime }\right) ^{3} &  & 
\end{array}%
\end{equation}%
and for $0<m\leq x\leq M$, 
\begin{equation}
\begin{array}{lcl}
\label{mge21}\mu _{3}\geq \frac{\mu _{2}-\mu _{1}^{\prime ^{2}}}{\mu
_{1}^{\prime }}\mu _{2}. &  & 
\end{array}%
\end{equation}
\end{theorem}

\begin{proof}
From the first inequality \eqref{ge3}, we get that 
\begin{equation}
\begin{array}{lcl}
\label{mge22}\mu _{2}\mu _{3}\geq \frac{\mu _{2}^{3}-\left( \mu _{1}^{\prime
}-m\right) ^{2}\mu _{2}^{2}}{\mu _{1}^{\prime }-m}. &  & 
\end{array}%
\end{equation}%
Let $f\left( x\right) =x^{3}-k^{2}x^{2}$. Then $f\left( x\right) $ has
minimum at $x=\frac{2}{3}k^{2}$ where its value is $-\frac{4}{27}k^{6}.$
Thus we have 
\begin{equation}
\begin{array}{lcl}
\label{mge23}\mu _{2}^{3}-\left( \mu _{1}^{\prime }-m\right) ^{2}\mu
_{2}^{2}\geq -\frac{4}{27}\left( \mu _{1}^{\prime }-m\right) ^{6}. &  & 
\end{array}%
\end{equation}%
Combining \eqref{mge22} and \eqref{mge23}; we immediately get the first
inequality \eqref{mge19}. The second inequality \eqref{mge19} follows
similarly from the second inequality \eqref{ge3}.

\noindent To prove \eqref{mge20} we write \eqref{ge3} in the form 
\begin{equation}
\begin{array}{lcl}
\label{mge24}m\mu _{2}^{\prime }+\frac{\left( m\mu _{1}^{\prime }-\mu
_{2}^{\prime }\right) ^{2}}{\mu _{1}^{\prime }-m}\leq \mu _{3}^{\prime }\leq
M\mu _{2}^{\prime }-\frac{\left( M\mu _{1}^{\prime }-\mu _{2}^{\prime
}\right) ^{2}}{M-\mu _{1}^{\prime }}. &  & 
\end{array}%
\end{equation}%
The inequalities \eqref{mge20} \ follow from \eqref{mge24} and the fact that
the right hand side and left hand side expressions \eqref{mge24} achieve
their maxima and minima at $\mu _{2}^{\prime }=\frac{M\left( M+\mu
_{1}^{\prime }\right) }{2}$ and $\mu _{2}^{\prime }=\frac{m\left( m+\mu
_{1}^{\prime }\right) }{2}$, respectively.

\noindent Likewise, the inequalities \eqref{mage4} follow from \eqref{ge3}
and the fact that the right hand side and left hand side expressions %
\eqref{ge3} achieve their maxima and minima at $\mu _{2}=\frac{\left( M-\mu
_{1}^{\prime }\right) ^{2}}{2}$ and $\mu _{2}=\frac{\left( \mu _{1}^{\prime
}-m\right) ^{2}}{2}$, respectively.

\noindent Further, from the first inequality \eqref{ge3} we have 
\begin{equation}
\begin{array}{lcl}
\label{mge26}m\leq \mu _{1}^{\prime }-\frac{\sqrt{\mu _{3}^{2}+4\mu _{2}^{3}}%
-\mu _{3}}{2\mu _{2}}. &  & 
\end{array}%
\end{equation}

\noindent For $m\geq 0,$ the inequalities \eqref{mge21} and \eqref{mge26}
are equivalent.
\end{proof}

\bigskip \noindent The inequalities \eqref{mge27} give the upper bounds for $%
m_{2}$ in terms of the values of $m,\ M\ $and $\overline{x}$. We prove
analogous inequalities involving $m_{2}$ and $m_{4}$.

\begin{theorem}
\label{rt6} For $m\leq x_{i}\leq M,\ i=1,2,...,n$, the inequalities 
\begin{equation}
\begin{array}{lcl}
\label{mge28}\frac{m_{2}^{4}}{m_{4}}\leq \frac{\left( n-1\right) ^{3}}{%
n^{2}-3n+3}\left( \overline{x}-m\right) ^{4} &  & 
\end{array}%
\end{equation}%
and 
\begin{equation}
\begin{array}{lcl}
\label{mge29}\frac{m_{2}^{4}}{m_{4}}\leq \frac{\left( n-1\right) ^{3}}{%
n^{2}-3n+3}\left( M-\overline{x}\right) ^{4} &  & 
\end{array}%
\end{equation}%
hold true.
\end{theorem}

\begin{proof}
For $r=2$, it follows from \eqref{ge5} that 
\begin{equation}
\begin{array}{lcl}
\label{mge30}\frac{(M-\overline{x})^{4}}{m_{4}}\leq \frac{\left( n-1\right)
^{3}}{n^{2}-3n+3}. &  & 
\end{array}%
\end{equation}%
From \eqref{age1}, we have 
\begin{equation}
\begin{array}{lcl}
\label{mge31}\frac{m_{2}}{M-\overline{x}}\leq \overline{x}-m. &  & 
\end{array}%
\end{equation}%
Combining \eqref{mge30} and \eqref{mge31} we immediately get \eqref{mge28}.
The inequality \eqref{mge29} follows on using similar arguments.
\end{proof}

\bigskip \noindent The inequalities \eqref{mge27} can equivalently be
written respectively in the form%
\begin{equation}
\begin{array}{lcl}
\label{mage7} m\leq \overline{x}-\frac{\sqrt{m_{2}}}{\sqrt{n-1}}\text{ \ and
\ }M\geq \overline{x}+\frac{\sqrt{m_{2}}}{\sqrt{n-1}}, &  & 
\end{array}%
\end{equation}%
and are important in finding the upper (lower) bound for the smallest
(largest) value of the data. We mention here analogous inequalities
involving $m_{2}$ and $m_{4}.$

\begin{corollary}
\label{cr3} Under the conditions of the above theorem, we have 
\begin{equation}
\begin{array}{lcl}
\label{mge32}m\leq \overline{x}-\left( \frac{n^{2}-3n+3}{\left( n-1\right)
^{3}}\frac{m_{2}^{4}}{m_{4}}\right) ^{\frac{1}{4}} &  & 
\end{array}%
\end{equation}%
and 
\begin{equation}
\begin{array}{lcl}
\label{mge33}M\geq \overline{x}+\left( \frac{n^{2}-3n+3}{\left( n-1\right)
^{3}}\frac{m_{2}^{4}}{m_{4}}\right) ^{\frac{1}{4}}. &  & 
\end{array}%
\end{equation}
\end{corollary}

\begin{proof}
The inequalities \eqref{mge32} and \eqref{mge33} follow immediately from the
inequalities \eqref{mge28} and \eqref{mge29}, respectively.
\end{proof}

\noindent The inequalities \eqref{mge32} and \eqref{mge33} give better
estimates than the corresponding estimates given by the inequalities %
\eqref{mage7}. Note that the inequality 
\begin{equation*}
\overline{x}-\left( \frac{n^{2}-3n+3}{\left( n-1\right) ^{3}}\frac{m_{2}^{4}%
}{m_{4}}\right) ^{\frac{1}{4}}\leq \overline{x}-\frac{\sqrt{m_{2}}}{\sqrt{n-1%
}}
\end{equation*}%
holds true if and only if%
\begin{equation*}
\frac{m_{4}}{m_{2}^{2}}\leq \frac{n^{2}-3n+3}{n-1}.
\end{equation*}%
This is true, see \cite{raj18}.

\noindent The well known Karl Pearson coefficient of dispersion $V=\frac{%
\sqrt{m_{2}}}{\overline{x}}$ is a widely used measure of dispersion. We
mention a lower for the ratio $\frac{\alpha _{4}}{V^{4}}$.

\begin{corollary}
\label{cr4} For $x_{i}>0,\ i=1,2,...,n$, we have 
\begin{equation}
\begin{array}{lcl}
\label{mge34}\frac{\alpha _{4}}{V^{4}}\geq \frac{n^{2}-3n+3}{\left(
n-1\right) ^{3}}. &  & 
\end{array}%
\end{equation}
\end{corollary}

\begin{proof}
If all the $x_{i}$'s are positive, $m\geq 0$, and \eqref{mge34} follows from %
\eqref{mge32}.
\end{proof}

\section{Bounds for eigenvalues and spreads of matrices}

\noindent Let $\mathbb{M}(n)$ denote the algebra of all $n\times n$ complex
matrices. The eigenvalues $\lambda \left( A\right) $ of an element $A\in 
\mathbb{M}(n)$ are the roots of the characteristic polynomial $\det \left(
A-\lambda \left( A\right) I\right) =0$ and are difficult to evaluate in
general. The bounds for eigenvalues have been studied extensively in
literature, for example see \cite{bha97, raj10, raj15, wol80}. In
particular, Wolkowicz and Styan \cite{wol80} have shown that%
\begin{equation}
\begin{array}{lcl}
\label{magen1}\lambda _{\min }(A)\leq \frac{\text{tr}A}{n}-\sqrt{\frac{\text{%
tr}B^{2}}{n(n-1)}}\text{ \ and \ }\lambda _{\max }(A)\geq \frac{\text{tr}A}{n%
}+\sqrt{\frac{\text{tr}B^{2}}{n(n-1)}}, &  & 
\end{array}%
\end{equation}%
where $B=A-\frac{\text{tr}A}{n}I$. \ We prove an extension of this result in
the following theorem.

\begin{theorem}
\label{rt7} If the eigenvalues of $A\in \mathbb{M}(n)$ are all real, then 
\begin{equation}
\begin{array}{lcl}
\label{mgen1}\lambda _{\min }(A)\leq \frac{\text{tr}A}{n}-\left( \frac{%
n^{2}-3n+3}{n^{3}\left( n-1\right) ^{3}}\right) ^{\frac{1}{4}}\frac{\text{tr}%
B^{2}}{\left( \text{tr}B^{4}\right) ^{\frac{1}{4}}} &  & 
\end{array}%
\end{equation}%
and%
\begin{equation}
\begin{array}{lcl}
\label{mgen2}\lambda _{\max }(A)\geq \frac{\text{tr}A}{n}+\left( \frac{%
n^{2}-3n+3}{n^{3}\left( n-1\right) ^{3}}\right) ^{\frac{1}{4}}\frac{\text{tr}%
B^{2}}{\left( \text{tr}B^{4}\right) ^{\frac{1}{4}}}. &  & 
\end{array}%
\end{equation}
\end{theorem}

\begin{proof}
The arithmetic mean of the eigenvalues $\lambda _{i}(A)$ can be written as $%
\overline{\lambda }\left( A\right) =\frac{1}{n}\sum_{i=1}^{n}\lambda _{i}(A)=%
\frac{\text{tr}A}{n}$. The second and fourth central moment of the
eigenvalues can be expressed in terms of tr$B^{2}$ and tr$B^{4}$,
respectively. We have 
\begin{equation}
\begin{array}{lcl}
\label{mgen3}m_{2}=\frac{1}{n}\sum_{i=1}^{n}\left( \lambda _{i}(A)-\overline{%
\lambda }\left( A\right) \right) ^{2}=\frac{\text{tr}B^{2}}{n}\text{ \ and \ 
}m_{4}=\frac{1}{n}\sum_{i=1}^{n}\left( \lambda _{i}(A)-\overline{\lambda }%
\left( A\right) \right) ^{4}=\frac{\text{tr}B^{4}}{n}. &  & 
\end{array}%
\end{equation}%
Apply Corollary 2.9., the inequalities \eqref{mgen1} and \eqref{mgen2}
follow on substituting the values of $m_{2}$ and $m_{4}$ from \eqref{mgen3}
and $\overline{x}=\overline{\lambda }\left( A\right) =\frac{\text{tr}A}{n}$
in \eqref{mge32} and \eqref{mge33}, respectively.
\end{proof}

\noindent The ratio of the largest and smallest eigenvalue of a positive
definite matrix $A$ is known as the ratio spread or condition number ($%
c\left( A\right) $) of a positive definite matrix. Wolkowicz and Styan \cite%
{wol80} have shown that%
\begin{equation}
\begin{array}{lcl}
\label{magen2}c\left( A\right) \geq 1+\frac{2\sqrt{\frac{\text{tr}B^{2}}{%
n(n-1)}}}{\frac{\text{tr}A}{n}-\sqrt{\frac{\text{tr}B^{2}}{n(n-1)}}}. &  & 
\end{array}%
\end{equation}%
It follows from Theorem $3.1$ that if $A$ is positive definite, then%
\begin{equation}
\begin{array}{lcl}
\label{mgen4}c\left( A\right) \geq 1+\frac{2}{\left( \frac{\left( n-1\right)
^{3}}{n\left( n^{2}-3n+3\right) }\right) ^{\frac{1}{4}}\frac{\left( \text{tr}%
B^{4}\right) ^{\frac{1}{4}}\text{tr}A}{\text{tr}B^{2}}-1}. &  & 
\end{array}%
\end{equation}

\begin{example}
Let%
\begin{equation*}
A_{1}=\left[ 
\begin{array}{cccc}
4 & 0 & 2 & 3 \\ 
0 & 5 & 0 & 1 \\ 
2 & 0 & 6 & 0 \\ 
3 & 1 & 0 & 7%
\end{array}%
\right] .
\end{equation*}%
From \eqref{magen1}, $\lambda _{\min }(A_{1})\leq 3.8417$ and $\lambda \max
(A_{1})\geq 7.1583$ while from \eqref{mgen1} and \eqref{mgen2} we have
better estimates $\lambda _{\min }(A_{1})\leq 3.7414$ and $\lambda \max
(A_{1})\geq 7.2586$, respectively. Also, from \eqref{magen2} and %
\eqref{mgen4}, we respectively have $c\left( A_{1}\right) \geq 1.8633$ and $%
c\left( A_{1}\right) \geq 1.9401$.
\end{example}

\begin{theorem}
\label{rt9} Under the conditions of Theorem $3.1$, we have 
\begin{equation}
\begin{array}{lcl}
\label{mgen11}\text{spd}(A)\geq \left( \frac{4}{n^{2}}\left( n\text{tr}B^{%
\text{4}}+3\left( \text{tr}B^{\text{2}}\right) ^{2}\right) \right) ^{\frac{1%
}{4}}. &  & 
\end{array}%
\end{equation}%
When $n$ is odd the inequality \eqref{mgen11} may be strengthened to 
\begin{equation}
\begin{array}{lcl}
\label{mgen12}\text{spd}(A)\geq \left( \frac{4}{n^{2}-1}\left( n\text{tr}B^{%
\text{4}}+3\left( \text{tr}B^{\text{2}}\right) ^{2}\right) \right) ^{\frac{1%
}{4}}. &  & 
\end{array}%
\end{equation}
\end{theorem}

\begin{proof}
The inequalities \eqref{mgen11} and \eqref{mgen12} follow respectively from
the inequalities \eqref{mge7} and \eqref{mge10} on using arguments similar
to those used in the proof of Theorem 3.1.
\end{proof}

\noindent Likewise, from the inequality \eqref{mge13}, we have 
\begin{equation}
\begin{array}{lcl}
\label{mgen13}\text{spd}(A)\geq \left( \frac{n^{3}}{\max_{j\in I}j\left(
n-j\right) \left( n^{2}-3nj+3j^{2}\right) }\text{tr}B^{\text{4}}\right) ^{%
\frac{1}{4}}, &  & 
\end{array}%
\end{equation}%
where $I=\left\{ 1,2,...,n-1\right\} $.

\noindent Also, from \eqref{mage1}, we have%
\begin{equation}
\begin{array}{lcl}
\label{magen3}\text{spd}(A)\geq 3\left( \frac{\text{tr}B^{\text{2}}\text{tr}%
B^{\text{4}}}{12n^{2}}\right) ^{\frac{1}{6}}. &  & 
\end{array}%
\end{equation}%
It is shown in \cite{raj18} that%
\begin{equation}
\begin{array}{lcl}
\label{magen4}\text{spd}(A)\geq \left( \frac{12\text{tr}B^{\text{4}}}{n}%
\right) ^{\frac{1}{4}}. &  & 
\end{array}%
\end{equation}
For $n=3$, \eqref{mgen12} and \eqref{mgen13} are identical. The inequality %
\eqref{mgen13} provides better estimate than \eqref{magen4}. The
inequalities \eqref{mgen11} and \eqref{mgen12} are clearly independent. We
show by means of examples that \eqref{mgen11}, \eqref{magen3} and %
\eqref{magen4} are independent.

\begin{example}
Let%
\begin{equation*}
A_{2}=\left[ 
\begin{array}{cccc}
1 & 1 & 0 & 2 \\ 
0 & 4 & 0 & 0 \\ 
0 & 3 & 1 & 1 \\ 
2 & 1 & 2 & 4%
\end{array}%
\right] .
\end{equation*}%
The eigenvalues of \ $A_{2}$ are $1,\ 4,\ \frac{5}{2}\pm \frac{\sqrt{33}}{2}$%
. From \eqref{mgen11} and \eqref{magen4}, spd$(A_{1})\geq 6.0264$, spd$%
(A_{2})\geq 4.6734$ and spd$(A_{1})\geq 6.2303$, spd$(A_{2})\geq 4.5767$,
respectively. So, the inequalities \eqref{mgen11} and \eqref{magen4} are
independent. From \eqref{magen3}, spd$(A_{1})\geq 6.3068$ and spd$%
(A_{2})\geq 4.762$. But for the matrix $A_{3}$ with eigenvalues $-1,\ 0,1$
with respective multiplicities $2,\ 5$ and $2$, we have \eqref{magen3} and %
\eqref{magen4} from spd$(A_{3})\geq 1.5131$ and spd$(A_{3})\geq 1.5197$,
respectively. From the inequality \eqref{mgen13}, we have spd$(A_{1})\geq
6.2549,$ spd$(A_{2})\geq 4.5948$ and spd$(A_{3})\geq 1.5902$. The inequality %
\eqref{mgen13} makes an improvement on \eqref{magen4}.
\end{example}

\noindent We now show that some of the above inequalities can be extended
for positive linear functionals. A linear functional $\varphi :\mathbb{M}%
(n)\rightarrow \mathbb{C}$ is called positive if $\varphi \left( A\right)
\geq 0$ whenever $A\geq O$ and unital if $\varphi \left( I\right) =1$, see 
\cite{bha07, raj18}.

\begin{theorem}
\label{rt8} Let $\varphi :\mathbb{M}(n)\rightarrow \mathbb{C}$ be a positive
unital linear functional and let $A$ be any Hermitian element of $\mathbb{M}%
(n)$. Then 
\begin{equation}
\begin{array}{lcl}
\label{mgen5}\text{spd}(A)^{2}\geq {2}\sqrt{\varphi \left( B^{4}\right)
+3\varphi \left( B^{2}\right) ^{2}}, &  & 
\end{array}%
\end{equation}%
where $B=A-\varphi \left( A\right) I$.
\end{theorem}

\begin{proof}
By the spectral theorem, for $r=1,2,...,$ we have 
\begin{equation}
\begin{array}{lcl}
\label{mgen6}A^{r}=\sum_{i=1}^{n}\lambda _{i}^{r}\left( A\right) P_{i}\text{
\ \ and \ \ }B^{r}=\sum_{i=1}^{n}\left( \lambda _{i}\left( A\right) -\varphi
\left( A\right) \right) ^{r}P_{i}, &  & 
\end{array}%
\end{equation}%
where $P_{i}$ are corresponding projections, $P_{i}\geq O$ and $%
\sum_{i=1}^{n}P_{i}=I$.

\noindent On applying $\varphi $, we find from \eqref{mgen6} that 
\begin{equation}
\begin{array}{lcl}
\label{mgen7}\varphi \left( A^{r}\right) =\sum_{i=1}^{n}\lambda
_{i}^{r}\left( A\right) \varphi \left( P_{i}\right) \text{ \ \ and \ \ }%
\varphi \left( B^{r}\right) =\sum_{i=1}^{n}\left( \lambda _{i}\left(
A\right) -\varphi \left( A\right) \right) ^{r}\varphi \left( P_{i}\right) & 
& 
\end{array}%
\end{equation}%
with $\sum_{i=1}^{n}\varphi \left( P_{i}\right) =1$.

\noindent Note that $\varphi \left( A\right) ,$ $\varphi \left( B^{2}\right) 
$ and $\varphi \left( B^{4}\right) $ are respectively the arithmetic mean,
second and fourth central moments of the eigenvalues $\lambda _{i}\left(
A\right) $ with respective weights $\varphi \left( P_{i}\right) ,\
i=1,2,...,n$. So, we can apply Corollary 2.2, and the inequality %
\eqref{mgen5} follows from \eqref{mge7}.
\end{proof}

\noindent On using similar arguments one can easily obtain from Theorem 2.5
that%
\begin{equation}
\begin{array}{lcl}
\label{mgen8}\text{spd}(A)^{2}\geq \left( \frac{3^{5}}{4}\varphi \left(
B^{2}\right) \varphi \left( B^{4}\right) \right) ^{\frac{1}{3}}. &  & 
\end{array}%
\end{equation}%
Further, it is shown in \cite{raj18} that 
\begin{equation*}
\begin{array}{lcl}
\label{mgen9}\text{spd}(A)^{3}\geq 6\sqrt{3}\varphi \left( B^{3}\right) , & 
& 
\end{array}%
\end{equation*}%
and by \eqref{mge16} we have the following refinement of this inequality%
\begin{equation}
\begin{array}{lcl}
\label{mgen10}\text{spd}(A)^{3}\geq 6\sqrt{3}\left( \varphi \left(
B^{2}\right) \varphi \left( B^{4}\right) -\varphi \left( B^{2}\right)
^{3}\right) ^{\frac{1}{2}}\geq 6\sqrt{3}\varphi \left( B^{3}\right) . &  & 
\end{array}%
\end{equation}

\section{Bounds for the roots and span of a polynomial}

\noindent Some bounds for the roots and span of a polynomial with all its
roots real are considered in \cite{raj18}. We here obtain some more bounds
for the roots and span in terms of the first five coefficients of the
polynomial.

\bigskip \noindent Let $x_{1},\ x_{2},...,x_{n}$ denote the roots of the
monic polynomial equation%
\begin{equation}
\begin{array}{lcl}
\label{pgen1}f(x)=x^{n}+a_{2}x^{n-2}+a_{3}x^{n-3}+...+a_{n-1}x+a_{n}=0. &  & 
\end{array}%
\end{equation}%
We assume that all the roots of $f(x)$ are real. On using the relations
between roots and coefficients of a polynomial, one can see that the
arithmetic mean $\overline{x}$ of the $x_{i}$'s equals zero. The second and
fourth central moment can respectively be written as, see \cite{raj18},%
\begin{equation}
\begin{array}{lcl}
\label{pgen2}m_{2}=-\frac{2}{n}a_{2}\text{ \ and \ }m_{4}=\frac{2}{n}\left(
a_{2}^{2}-2a_{4}\right) . &  & 
\end{array}%
\end{equation}

\begin{theorem}
Let the roots of the polynomial \eqref{pgen1} be all real. Denote by $%
x_{\min }(f)$ and $x_{\max }(f)$ the smallest and largest root of $f(x)=0$.
Then, for $n\geq 5,$ we have%
\begin{equation}
\begin{array}{lcl}
\label{pgen3}x_{\min }(f)\leq -\left( \frac{8\left( n^{2}-3n+3\right) }{%
n^{3}\left( n-1\right) ^{3}}\frac{a_{2}^{4}}{a_{2}^{2}-2a_{4}}\right) ^{%
\frac{1}{4}} &  & 
\end{array}%
\end{equation}%
and%
\begin{equation}
\begin{array}{lcl}
\label{pgen4}x_{\max }(f)\geq \left( \frac{8\left( n^{2}-3n+3\right) }{%
n^{3}\left( n-1\right) ^{3}}\frac{a_{2}^{4}}{a_{2}^{2}-2a_{4}}\right) ^{%
\frac{1}{4}}. &  & 
\end{array}%
\end{equation}
\end{theorem}

\begin{proof}
Apply Corollary 2.9., the inequalities \eqref{pgen3} and \eqref{pgen4}
follow on substituting the values of $m_{2}$ and $m_{4}$ from \eqref{pgen2}
and $\overline{x}=0$ in \eqref{mge32} and \eqref{mge33}, respectively.
\end{proof}

\begin{theorem}
If roots of the polynomial \eqref{pgen1} are all real, then%
\begin{equation}
\begin{array}{lcl}
\label{pgen5}\text{spn}\left( f\right) \geq \left( 8\left( \frac{%
a_{2}^{2}-2a_{4}}{n}+\frac{6a_{2}^{2}}{n^{2}}\right) \right) ^{\frac{1}{4}},
&  & 
\end{array}%
\end{equation}%
\begin{equation}
\begin{array}{lcl}
\label{pgen6}\text{spn}\left( f\right) \geq \left( \frac{%
2n^{3}(a_{2}^{2}-2a_{4})}{\max_{j\in I}j\left( n-j\right) \left(
n^{2}-3nj+3j^{2}\right) }\right) ^{\frac{1}{4}} &  & 
\end{array}%
\end{equation}%
and 
\begin{equation}
\begin{array}{lcl}
\label{pgen7}\text{spn}\left( f\right) \geq \left( \frac{3^{5}}{n^{2}}%
a_{2}\left( 2a_{4}-a_{2}^{2}\right) \right) ^{\frac{1}{6}}. &  & 
\end{array}%
\end{equation}%
When $n$ is odd the inequality \eqref{pgen5} may be strengthened to%
\begin{equation}
\begin{array}{lcl}
\label{pgen8}\text{spn}\left( f\right) \geq \left( \frac{8n}{n^{2}-1}\left(
a_{2}^{2}-2a_{4}+\frac{6}{n}a_{2}^{2}\right) \right) ^{\frac{1}{4}}. &  & 
\end{array}%
\end{equation}
\end{theorem}

\begin{proof}
On using arguments similar to those used in the proof of Theorem 4.1 the
inequalities \eqref{pgen5}, \eqref{pgen6}, \eqref{pgen7} and \eqref{pgen8}
follow respectively from \eqref{mge7}, \eqref{mge13}, \eqref{mage1} and %
\eqref{mge10}.
\end{proof}

\noindent \vskip0.2in \noindent \textbf{Acknowledgments} The authors are
grateful to Prof. Rajendra Bhatia for the useful discussions and
suggestions, and first author thanks Ashoka University for a visit in
January 2019. The support of UGC-SAP is acknowledged.


\end{document}